\documentclass[transmag]{IEEEtran}
\usepackage{latexsym}
\usepackage{graphicx}
\usepackage{amsfonts,amssymb,amsmath}
\usepackage{hyperref}

\usepackage{scalerel}[2016/12/29]

\def\BibTeX{{\rm B\kern-.05em{\sc i\kern-.025em b}\kern-.08em T\kern-.1667em\lower.7ex\hbox{E}\kern-.125emX}}

\begin{document}

\bibliographystyle{IEEEtran}

\newcommand{\dl}{\delta}
\newcommand{\reals}{\mathbb{R}}
\newcommand{\vphi}{\varphi}
\newcommand{\F}{\mathcal{F}}
\newcommand{\HH}{\mathcal{H}}
\newcommand{\pullback}{{\scaleobj{0.65}{<}}}

\newcommand{\lra}{\longrightarrow}

\title{A Note on the Convolution of Circle Impulses}

\author{Brad Osgood,
\IEEEmembership{Member, IEEE}
\thanks{Brad Osgood is a Professor of Electrical Engineering at Stanford University, Stanford, CA 94305 USA (e-mail: osgood@stanford.edu).}
}

\IEEEtitleabstractindextext{\begin{abstract}
We present a formula for the convolution $\dl_{C_1}*\dl_{C_2}$ of two circle impulses. The derivation uses a classical addition theorem for Bessel functions. 
\end{abstract}

\begin{IEEEkeywords}
Convolution, Impulses, Fourier transform,  Hankel transform,
\end{IEEEkeywords}
}
\maketitle

\section{INTRODUCTION and MAIN RESULT}

\IEEEPARstart{T}{he} purpose of this short note is to derive a formula for the convolution of two circle impulses. To my knowledge such a result has not appeared. Appealing aspects of the derivation are the use of an old addition formula for Bessel functions, and at one point the algebra is rather remarkable. We proceed with the formula.

For $i=1,2$ let $C_i$ be the circle $\|x-b_i\|=R_i$ centered at $b_i\in\reals^2$ with radius $R_i$. Let $\dl_{C_i}$ be the impulse supported on $C_i$ (definition below) and let $\rho=\| x-(b_1+b_2)\|$. The main result is that convolution results in the function
\begin{equation}\label{eq:dl*dl}
\small{(\dl_{C_1}*\dl_{C_2})(\rho)=
{\frac{4R_1R_2}{\sqrt{(\rho^2-(R_1-R_2)^2)((R_1+R_2)^2-\rho^2))}}}},
\end{equation}
if $|R_1-R_2| < \rho < R_1+R_2$, and $(\dl_{C_1}*\dl_{C_2})(\rho)
 =0$ if $0< \rho<|R_1-R_2|$ or $\rho>R_1+R_2$.
In particular, the formula expresses the fact that $\dl_{C_1}*\dl_{C_2}$ is radial with respect to the point $b_1+b_2$. It is also symmetric in $C_1$ and $C_2$ as it should be. 
 
\section{DEFINITIONS and BACKGROUND}

One can define a circle impulse $\dl_C$ in several equivalent ways: via a limiting process, as for example in \cite{bracewell:imaging}, or directly as a distribution paired with a smooth test function $\vphi$ that is of compact support on $\reals^2$. The pairing is
\begin{equation} \label{eq:dl_C_paired_vphi}
 \langle \dl_C,\vphi \rangle = \int_{C} \vphi\,ds,
 \end{equation}
where $ds$ is arclength along $C$.

\subsection{Use of Pullbacks}

It is useful for definitions and formulas to employ pullbacks of functions and distributions, see \cite{friedlander:distributions} for the general setup and \cite{osgood:fourier} for special cases. For the circle $\|x-b\|=R$, the formula \eqref{eq:dl_C_paired_vphi} in fact comes from the definition
of $\dl_C$ as the pullback, $\dl_C = g^\pullback \dl$, of the usual one-dimensional point $\dl$ by $g:\reals^2\lra \reals$, $g(x) = \|x-b\|-R$. Keeping to the custom of evaluating $\dl$'s at points, this is usually expressed as $\dl_C(x_1,x_2)=\dl(g(x_1,x_2))$. Note that $\dl_C$ and subsequent formulas depend on the choice of the function $g$  describing $C$. 
Observe that for pullback we write $g^\pullback$ instead of the more common $g^*$ so as not to clash,  inevitably, with $*$ used for convolution.

Let $\tau_b(x)=x-b$ be the shift operator on $\reals^2$. It is easy to show  that if $C$ is the circle $\|x\|=R$ then $\tau_b^\pullback\dl_C$ is the circle impulse for $\|x-b\|=R$. Moreover, 
 as a shift of the convolution is the convolution of the shifts,   if $C_i$, $i=1,2$, are the circles $\|x\|=R_i$ we have
\begin{equation} \label{eq:shifted-convolution}
(\tau_{b_1+b_2})^\pullback(\dl_{C_1}*\dl_{C_2}) = \tau_{b_1}^\pullback \dl_{C_1}*\tau_{b_2}^\pullback\dl_{C_2},
\end{equation}
so it suffices to prove \eqref{eq:dl*dl} for circles concentric to the origin. 

\subsection{Radial Functions}

To express a function (distribution) $f$ in polar coordinates is to form $\tilde f(r,\theta) = (P^\pullback\! f)(r,\theta)$ using the polar coordinate mapping $P(r,\theta) = (r\cos\theta,r\sin\theta)$. A function (distribution) is radial if $A^\pullback f = f$ for any orthogonal transformation $A$, in which case we write the polar coordinate version simply as $\tilde f(r)$. The convolution of two radial functions (distributions) is radial and $\dl_C$ is radial for a circle centered at the origin. The Fourier transform of a radial function is also radial and 
\[
(\F f)^\sim(r) = \HH\tilde{f}(r) =  2\pi \int_0^\infty \tilde{f}(\rho) J_0(2\pi r \rho)\,\rho d\rho.
\] 
Here, $J_0$ is the zeroth order Bessel function of the first kind and $\HH$ is the Hankel transform.

For the circle $C$, $\| x\|=R$, the Fourier transform of $\dl_C$ is the radial function
$
(\F \dl_C)^\sim (r) = 2\pi R J_0(2\pi rR)
$.
This is a standard result. We can write it as 
\begin{equation*}
\HH \dl_C^\sim(r) = 2\pi R J_0(2\pi r R)
\end{equation*}

 See  \cite{grafakos:radial} for a general discussion of radial distributions. 

\section{CONVOLUTION and the NEUMANN ADDITION FORMULA}

We consider the convolution $\dl_{C_1}*\dl_{C_2}$ for the concentric circles $\|x\|=R_i$. Invoking the convolution theorem and then expressing the result in polar coordinates, we have $(\F(\dl_{C_1} * \dl_{C_2}))^\sim = (\F\dl_{C_1})^\sim(\F\dl_{C_2})^\sim$, or
\begin{equation} \label{eq:Hankel-product}
\begin{aligned}
\HH(\dl_{C_1} &* \dl_{C_2})^\sim(r) = \HH \dl^{\sim}_{C_1}(r)\HH\dl^{\sim}_{C_2}(r) \\
&= (2\pi)^2R_1R_2J_0(2\pi R_1r)J_0(2\pi R_2r).
\end{aligned}
\end{equation}
The Hankel transform is its own inverse, so we want to find $\HH$ applied to the right-hand side. 

An old result on the product of two Bessel functions, a consequence of the Neumann addition formula, comes to our aid. In the present setting it reads
\begin{equation} \label{eq:Neumann-1}
\begin{aligned}
J_0(2\pi R_1r)J_0(2\pi R_2 r) &=\\
\frac{1}{2\pi} \int_0^{2\pi} J_0((R_1^2+R_2^2-2R_1R_2& \cos\theta)^{1/2}2\pi r)\,d\theta.
\end{aligned}
\end{equation}
See \cite{watson:bessel}, Section 11.2.
Temporarily drop the factor $(2\pi)^2R_1R_2$ in \eqref{eq:Hankel-product}, and to further save space let
\[
\Psi(\theta)   = (R_1^2+R_2^2-2R_1R_2 \cos\theta)^{1/2}.
\]
Thus,
\[
\begin{aligned}
\HH\left(\frac{1}{2\pi}\int_0^{2\pi} J_0(\Psi(\theta)2\pi r)\,d\theta \right)(\rho) &=\\
\int_0^\infty\!\!\int_0^{2\pi} J_0(2\pi r \rho)&J_0(\Psi(\theta)2\pi r)\,r d\theta dr
\end{aligned} 
\]
Interchange the order of integration and for the inner integral appeal to the ``closure equation:"
\[
\int_0^\infty J_0(kr)J_0(k'r)\,r dr = \frac{1}{k}\dl(k-k'), \quad k,k' >0.
\] 
(The classical version of this, not in terms of $\dl$, goes back to 
Hankel, see \cite{watson:bessel}, Section 14.4.) Then,
\[
\begin{aligned}
\int_0^\infty J_0(2\pi \rho r)J_0(2\pi \Psi(\theta)r)\,r dr &= \frac{1}{2\pi \rho}\dl(2\pi \rho-2\pi\Psi(\theta))\\
&= \frac{1}{(2\pi)^2\rho} \dl(\rho - \Psi(\theta)).
\end{aligned}
\]
Finally, put the factor $(2\pi)^2R_1R_2$ back to obtain
\begin{equation} \label{eq:formula-integral}
(\dl_{C_1}*\dl_{C_2})^\sim(\rho)= \frac{R_1R_2}{\rho} \int_0^{2\pi}  \dl(\rho-\Psi(\theta))\,d\theta.
\end{equation}

\subsection{$\dl$ and Zeros}
We can further analyze the right-hand side of \eqref{eq:formula-integral} using the well known formula 
\begin{equation} \label{eq:dl(f(theta))}
\dl(\Phi(\theta)) = \sum_n \frac{\dl(\theta -\theta_n)}{|\Phi'(\theta_n)|},
\end{equation}
where 
\[
\Phi(\theta) = \rho-\Psi(\theta)
\]
 and the $\theta_n$ are the zeros of $\Phi(\theta)$. This may be derived via pullbacks, or by other means, see, e.g., \cite{osgood:fourier}. The formula holds so long as  $\Phi'(\theta_n) \ne 0$. The derivative is 
\begin{equation} \label{eq:dervative_of_f}
\Phi'(\theta) = -\frac{R_1R_2 \sin\theta}{(R_1^2+R_2^2-2R_1R_2\cos \theta)^{1/2}}.
\end{equation}
This is zero at $\theta=0,\pi$ and $2\pi$, and only there, and we will exclude these points by restricting $\rho$. 

The minimum of $\Psi(\theta)=(R_1^2+R_2^2 -2R_1R_2\cos\theta)^{1/2}$  occurs for $\theta = 0$ and $\theta = 2\pi$ and has the value $|R_1-R_2|$.  The maximum occurs at $\theta = \pi$ and has the value $R_1+R_2$. Thus if we restrict $\rho$ to 
\begin{equation} \label{eq:range_of_rho}
|R_1-R_2| < \rho < R_1+R_2,
\end{equation}
then within this range $\Phi(\theta)$ has two zeros with nonzero derivatives, say at $0 < \theta _1 < \pi$ and, $\pi < \theta_2 < 2\pi$. Evidently $\theta_2=2\pi -\theta_1$ and  the derivatives have equal absolute value. Strictly outside this range of $\rho$, either  $0<\rho <|R_1-R_2|$ or $\rho > R_1+R_2$, there are no zeros and $\dl(\Phi(\theta))=0$.

Now, $\Phi(\theta)=0$ when 
\begin{equation} \label{eq:cosine}
\cos\theta =\frac{R_1^2+R_2^2-\rho^2}{2R_1R_2}.
\end{equation}
At such a value, 
\[
\Phi'(\theta)^2 = \left(\frac{R_1R_2}{\rho}\right)^2 \sin^2\theta =\left(\frac{R_1R_2}{\rho}\right)^2(1-\cos^2\theta),
\]
and we can use \eqref{eq:cosine} to write
\[
\Phi'(\theta)^2 = \left(\frac{R_1R_2}{\rho}\right)^2\left(1- \left( \frac{R_1^2+R_2^2-\rho^2}{2R_1R_2}\right)^2\right).
\]
 This simplifies (!) to
 \[
\Phi'(\theta)^2 = \frac{1}{4\rho^2}(\rho^2-(R_1-R_2)^2)((R_1+R_2)^2-\rho^2).
\]
Note how directly the condition \eqref{eq:range_of_rho} comes into play.

Going back to \eqref{eq:dl(f(theta))}, again denoting the two zeros by $\theta_1$ and $\theta_2=2\pi - \theta _1$ we have
$|\Phi'(\theta_1)| = |\Phi'(\theta_2)|$ and
\[
\begin{aligned}
\dl(\Phi(&\theta))=
 \frac{2\rho(\dl(\theta-\theta_1)+\dl(\theta-\theta_2))}{\sqrt{(\rho^2-(R_1-R_2)^2)((R_1+R_2)^2-\rho^2))}},
 \end{aligned}
\]
from which
\[
\begin{aligned}
\frac{R_1R_2}{\rho} \int_0^{2\pi}  \dl(\Phi(\theta))\,d\theta &= \\
&\hspace{-.2in}\frac{4R_1R_2}{\sqrt{(\rho^2-(R_1-R_2)^2)((R_1+R_2)^2-\rho^2))}}.
\end{aligned}
\]
This is \eqref{eq:dl*dl} for circles concentric at the origin. The general result follows from this as in \eqref{eq:shifted-convolution}

\subsection{Sample Plots}

Here is a surface plot of $\dl_{C_1}*\dl_{C_2}$ for the concentric circles $\|x\|=2$ (red) and $\|x\|= 3$ (blue), together with the profile plot.

\begin{center}
\includegraphics[scale=.4]{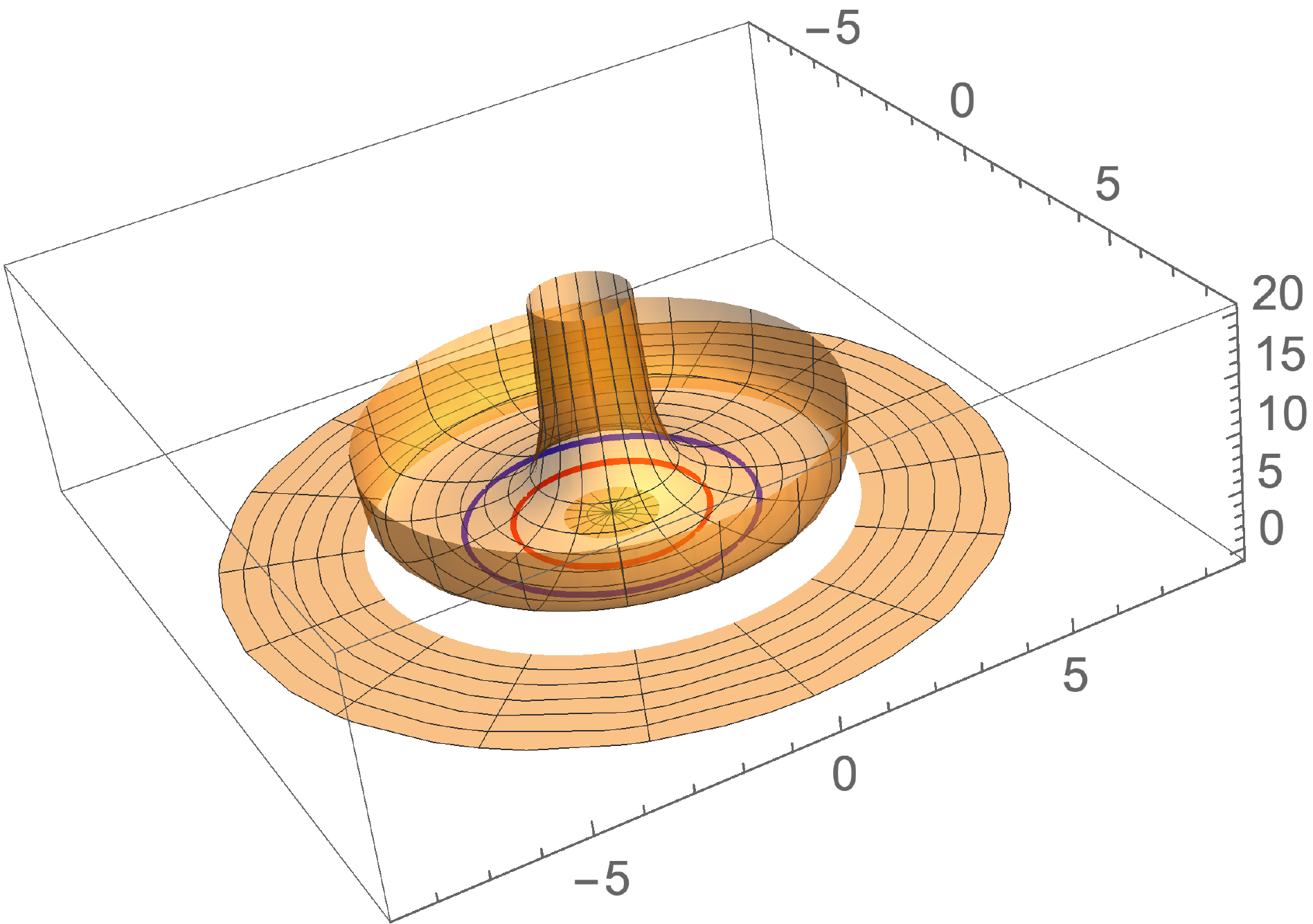}

\includegraphics[scale=.4]{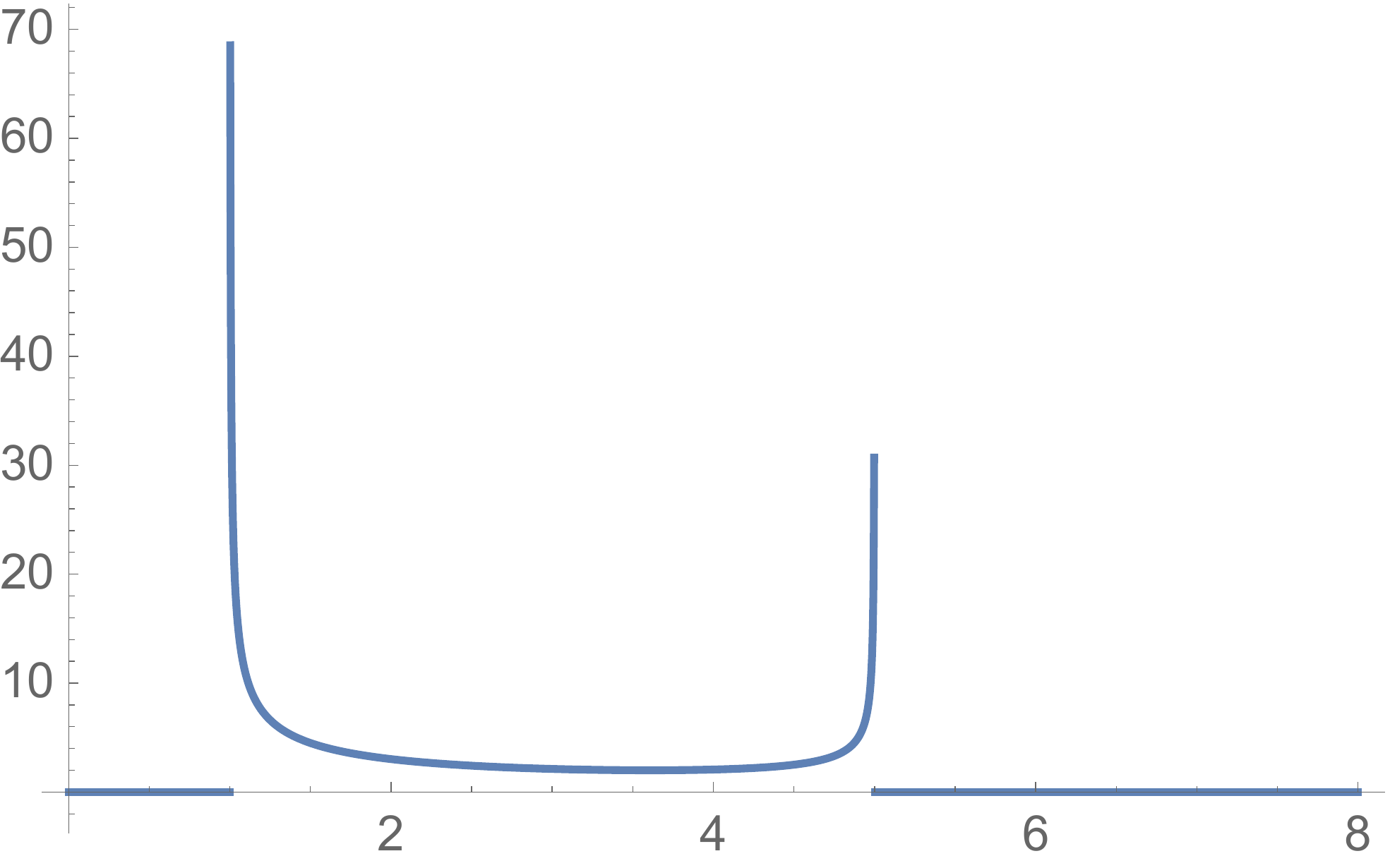}
\end{center}

It is  interesting to experiment with different configurations of circles.

\section{Additional Formulas}

Finally, it is natural to inquire about two other operations with impulses, namely multiplication and convolution with a function. Let $C$ be the circle $\|x\|=R$. For a function $f$ and for $x \ne 0$ let $f_C(x)=f(Rx/\|x\|)$. Then 
\begin{equation} \label{eq:f-times-dl}
f\dl_C = f_C \dl_C.
\end{equation}
Note that if $f$ is radial then $f\dl_C = \tilde{f}(R)\dl_C$.

For convolution the result is 
\begin{equation} \label{eq:f*dl}
(f*\dl_C)(x_1,x_2) = \int_0^{2\pi} f(x_1+R\cos\theta, x_2+R\sin\theta)\,Rd\theta.
\end{equation}
In words, forming $(f*\dl_C)(x)$ replaces $f(x)$ with an average of $f$ on the circle centered at $x$ of radius $R$. Thus $(f*\dl_C)(x_0)$ depends only on the values of $f$ on the circle $\|x-x_0\|=R$, a kind of ``locally radializing'' the function at each point.

The derivations of \eqref{eq:f-times-dl} and \eqref{eq:f*dl} are straightforward and we omit the details.


 \bibliography{references}


\end{document}